\documentclass[a4paper,twoside,11pt]{article}
\usepackage[english]{babel}
\usepackage[T1]{fontenc}
\usepackage[ansinew]{inputenc}
\usepackage{geometry}
\usepackage{color} 
\usepackage{lmodern}

\usepackage{graphicx}

\usepackage{amsmath} \numberwithin{equation}{section}
\usepackage{amsthm}
\usepackage{amsfonts}
\usepackage{amssymb}
\usepackage{graphics}
\usepackage{color}

\theoremstyle{plain}
\newtheorem{theorem}{Theorem}[section]

\newtheorem{remark}[theorem]{Remark}
\theoremstyle{definition}

\theoremstyle{remark}

\definecolor{orange}{rgb}{1,0.5,0}

\usepackage{enumerate}
\newcommand\C{\mathbb C}         
\newcommand\R{\mathbb R}

\newcommand\D{\mathbb D}

\renewcommand\Re{\text{Re }}

\newcommand{\arctanh}{\operatorname{arctanh}}

\usepackage[pdfstartview=FitH]{hyperref}

\begin{document}
\parindent 0pt 

\setcounter{section}{0}

\title{Value Ranges of Univalent Self-Mappings of the Unit Disc}
\author{Julia Koch  \and Sebastian Schlei\ss inger
       \thanks{Supported by the ERC grant  ``HEVO - Holomorphic Evolution Equations'' n. 277691.}}
\date{\today}
\maketitle
\abstract{We describe the value set $\{f(z_0)\,:\, f:\mathbb{D}\to\mathbb{D} \text{ univalent}, f(0)=0, f'(0)= e^{-T} \},$ where $\mathbb{D}$ denotes the unit disc and $z_0\in\mathbb{D}\setminus\{0\}$, $T>0,$ by applying  Pontryagin's maximum principle to the radial Loewner equation.} \\

{\bf Keywords:}  Pontryagin's maximum principle, radial Loewner equation, univalent functions\\
{\bf 2010 Mathematics Subject Classification:} 30C55, 30C80.

%\tableofcontents
\parindent 0pt

\section{Introduction and main result}

Given a bounded univalent function $f$ on a simply connected domain $\Omega\subsetneq\C$ and two distinct points $a, b\in\Omega$, it is quite natural to ask the question which values $f(b)$ can take if $f(a)$ and $f'(a)$ are prescribed.
Since the Riemann mapping theorem tells us that any such domain $\Omega$ can be mapped conformally onto the unit disc $\D=\{z\in\C \, : \, |z|<1\}$  such that $a$ is mapped to $0$, the problem can be restricted to the case of $\Omega=\D$ and $a=0$. \\

By multiplying with a real constant $\leq1$ and applying an automorphism of $\D$, we may assume $f:\D\to\D$ and $f(0)=0.$ Then the Schwarz lemma tells us that $|f'(0)|\leq 1$ and $|f'(0)|=1$ if and only if $f$ is the rotation $f(z)=f'(0)z.$ In order to describe the non-trivial case $|f'(0)|<1,$ we can restrict ourselves to the case $f'(0)\in(0,1)$ because of rotational symmetry. Thus we consider the set 
$$\mathcal S_T:=\{f:\D\to\D \text{ univalent}, f(0)=0, f'(0)= e^{-T}\}, \quad T>0.$$
 In this note, we will determine the value set 
$$V_T(z_0)= \{f(z_0)\,:\,  f\in\mathcal S_T\}, \quad z_0 \in \D\setminus\{0\}.$$ 

Variations of the set $V_T(z_0)= \{f(z_0)\,:\,  f\in\mathcal S_T\}$ have been determined by various authors, from the classical setting of the Schwarz  and Rogosinksi's lemma \cite{rogo1}, which concern itself with holomorphic functions $f:\D\to\D$, $f(0)=0$ that fulfil no further conditions, to a recent paper by Roth and Schlei{\ss}inger \cite{MR3262210} that determines the set $\mathcal V(z_0)= \{f(z_0)\,:\,  f\in\mathcal S\}$, with the class $\mathcal S:=\{f:\D\to\D \text{ univalent}, f(0)=0, f'(0)>0\}$. Note that $\mathcal V(z_0)=\cup_{T>0} V_T(z_0)$.\\

Our results are analogous to the results of Prokhorov and Samsonova \cite{MR3334955}, who study univalent self-mappings of the upper half-plane having the so called hydrodynamical normalization at the boundary point $\infty.$ Finally we note that in \cite{MR0453993}, the authors consider the set $\{\log(f(z_0)/z_0) \,:\, f:\D\to\C \text{ univalent},\phantom{i} f(0)=0, \phantom{i}|f(z)|\leq M\}$ for $M>0.$ We use a different and more straightforward approach to directly determine the set $V_T(z_0)$ by applying Pontryagin's maximum principle to the radial Loewner equation.
\\

In the following, for the sake of simplicity, we assume that $z_0\in(0,1)$; for other values of $z_0$, we just consider the function $z\mapsto e^{i\arg z_0}f\left(e^{-i\arg z_0} z\right)$ instead of $f$. \\

\begin{theorem}\label{thm:1}Let $z_0\in (0,1).$
For $x_0\in[-1,1]$ and $T>0$, let $r=r(T, x_0)$
be the (unique) solution to the equation
\begin{align*}
(1+x_0)(1-z_0)^2\log(1-r)+(1-x_0)(1+z_0)^2\log(1+r)-(1-2x_0z_0+z_0^2)\log r=\\
(1+x_0)(1-z_0)^2\log(1-z_0)+(1-x_0)(1+z_0)^2\log(1+z_0)-(1-2x_0z_0+z_0^2)\log e^{-T
}z_0
\end{align*}
and let 
$$ \sigma(T, x_0)=\frac{2(1-z_0^2)\sqrt{1-x_0^2}}{1-2x_0z_0+z_0^2}\left(\arctanh z_0 - \arctanh r(T, x_0)\right). $$
 Furthermore, for fixed $T\geq 0,$ define the two curves $C_+(z_0)$ and $C_-(z_0)$ by
\begin{align*}
C_{\pm}(z_0):=\left\{w_{\pm}(x_0):= r(T, x_0)e^{\pm i\sigma(T, x_0)}\,:\,  x_0\in[-1,1]\right\}.
\end{align*}

Then, if $\arctanh z_0<\frac{\pi}{2}$,  $V_T(z_0)$ is the closed region whose boundary consists of the two curves $C_+(z_0)$ and $C_-(z_0)$, which only intersect at $x_0\in\{-1,1\}.$\\

For $\arctanh z_0\geq\frac{\pi}{2}$, there are two different cases: First assume that $T$ is large enough that the equation
\begin{equation}\label{split}
\frac{2(1-z_0^2)\sqrt{1-x^2}}{1+2x z_0+z_0^2}\left(\arctanh z_0-\arctanh r(T,x)\right)=\pi
\end{equation}
admits a solution $x\in[-1,1]$. Then the curves  $C_+(z_0)$ and $C_-(z_0)$ intersect more than twice.
There is a $\chi\in(-1,1)$ such that $\widetilde C_{+}(z_0)\cup  \widetilde C_{-}(z_0)$ is a closed Jordan curve, where 
\begin{align*}
 \widetilde C_{\pm}(z_0)&:=\left\{w_{\pm}(x_0) \,:\,  x_0\in[\chi,1]\right\},
\end{align*}
and an $\aleph\in(-1,1)$ such that  $\widehat C_{+}(z_0)\cup  \widehat C_{-}(z_0)$ is a closed Jordan curve, where 
\begin{align*}
 \widehat C_{\pm}(z_0)&:=\left\{w_{\pm}(x_0) \,:\,   x_0\in[-1,\aleph]\right\}.
\end{align*}
Then $V_T(z_0)$ is the closed region whose boundary is $\widetilde C_{+}(z_0)\cup  \widetilde C_{-}(z_0) \cup \widehat C_{+}(z_0)\cup  \widehat C_{-}(z_0)$.\\ 
For smaller $T$ that do not admit a solution to (\ref{split}), the set $V_T(z_0)$ can be described exactly as in the case of $\arctanh z_0<\frac{\pi}{2}$.

\end{theorem}

The following figures show the evolution of the sets $V_T(z_0)$ over time. Note that $\arctanh z_0=\frac{\pi}{2} \iff z_0=\tanh(\pi/2)\approx0.917.$

\begin{figure}[h]
\begin{minipage}{0.5\textwidth}
\includegraphics[width=0.9\textwidth]{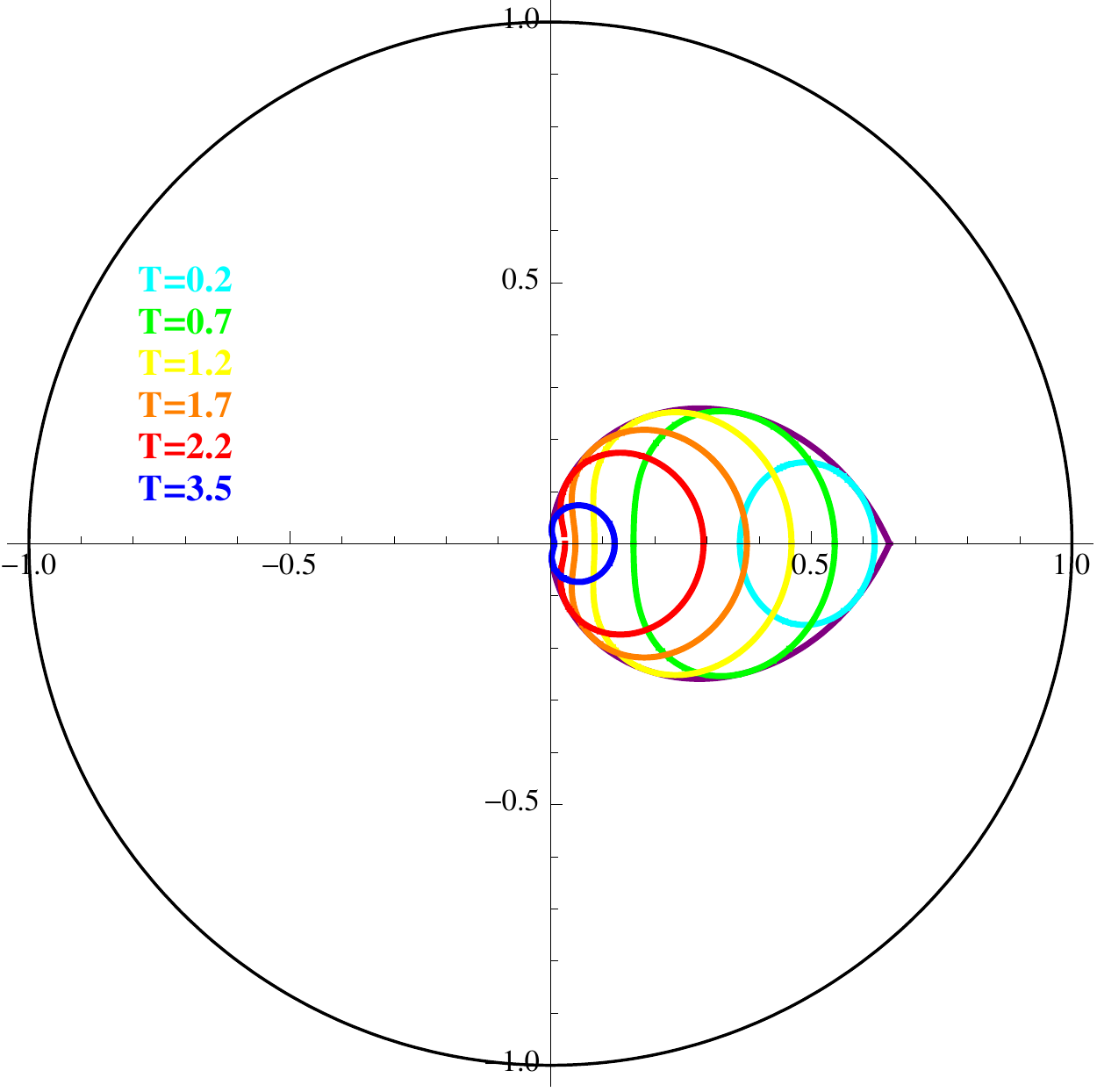}
\caption{$V_T(0.65)$}
\end{minipage}
\hfill
\begin{minipage}{0.5\textwidth}
\includegraphics[width=0.9\textwidth]{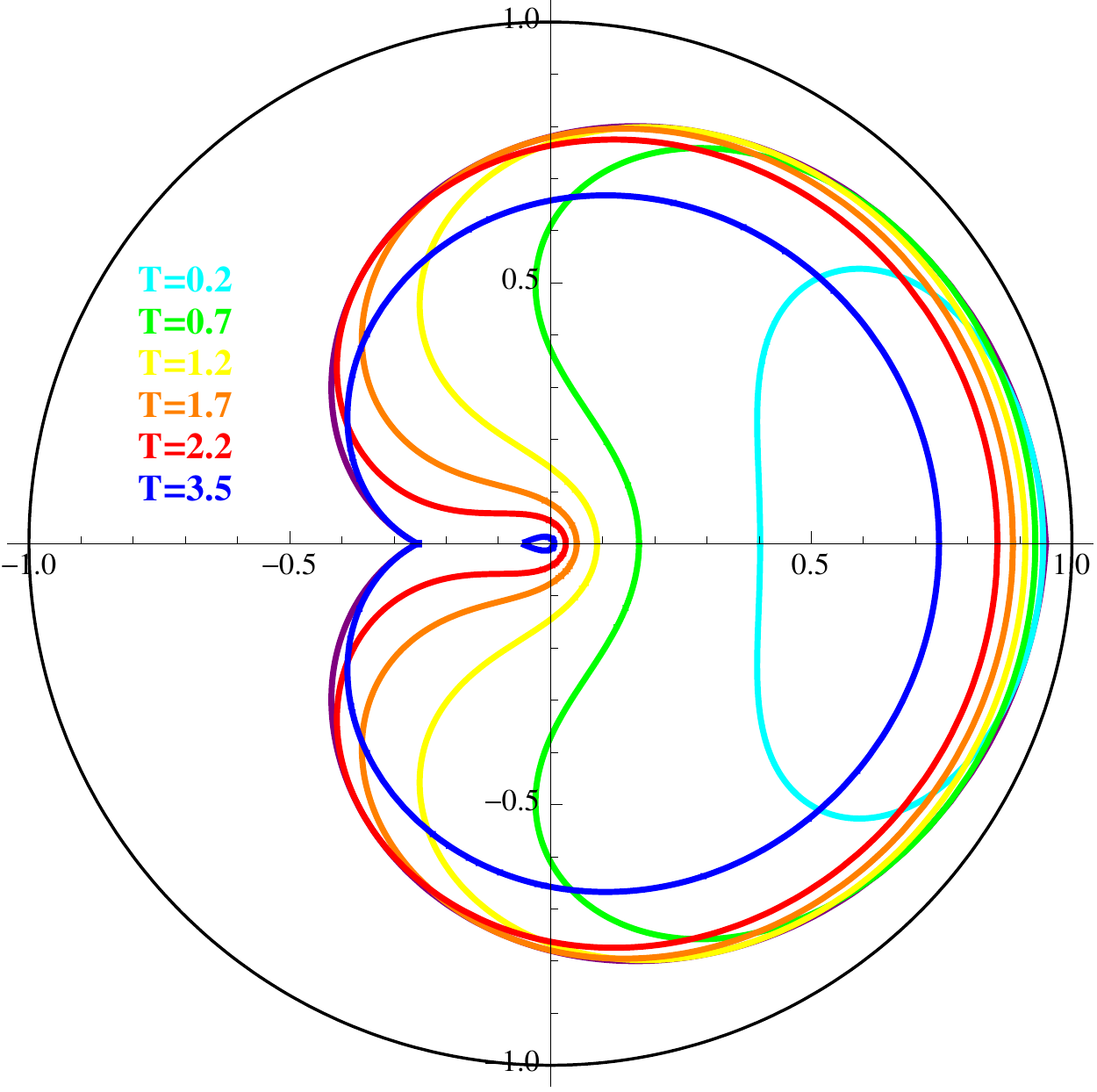}
\caption{$V_T(0.95)$}
\end{minipage}
\begin{center}
The sets $V_T(z_0)$ for $z_0=0.65, 0.95$ and 
$T=0.2+0.5j$, $j=0,1,\dots, 4$, and $T=3.5$.
\end{center}
\end{figure}

We prove Theorem \ref{thm:1} in Section \ref{sec:proof}, and in Section \ref{sec:inverse} we consider the similar problem of describing the value set $\{f^{-1}(z_0)\,:\, f\in \mathcal{S}_T \text{ with $z_0\in f(\D)$}\}$ for the inverse functions.

\section{Proof of Theorem \ref{thm:1}}\label{sec:proof}

Consider the radial Loewner equation
\begin{equation}\label{Loewner0}
\text{$\dot{f}_t(z) = -f_t(z) \cdot p(t, f_t(z))$ for a.e. $t\geq0$, \quad $f_0(z)=z\in\D$,} 
\end{equation}
%\begin{proof}
where $p:[0,\infty)\times \D\to \C$ is a \emph{Herglotz function}, i.e. for almost every $t\geq0$, $z\mapsto p(t,z)$ is a holomorphic function with $\Re p(z)>0$ for all $z\in\D$ and $p(0)=1$ and the function $t\mapsto p(t,z)$ is measurable for every $z\in\D.$ \\
For every $f\in \mathcal S_T$ there exists a Herglotz function $p(t,z)$ such that the solution $\{f_t\}_{t\geq0}$ of \eqref{Loewner0} satisfies $f_T=f$; see \cite{Pom:1975}, Chapter 6.\\

Thus the description of $V_T(z_0)$ can be translated into the control theoretic problem of describing the reachable set $R_T(z_0)$ of the initial value problem  
\begin{equation}\label{Loewner1} \dot{w}(t) = -w(t) \cdot p(t, w(t)),\quad w(0)=z_0\in\D, \end{equation}
where $p(t,z)$ runs through the set of all Herglotz functions and 
$$R_T(z_0):=\{w(T)\,:\, w:[0,T]\to\D \text{ solves \eqref{Loewner1}}\}.$$ 
Then we have $V_T(z_0) = R_T(z_0)$ and, obviously, $R_T(z_0)$ is a closed set.\\

Denote by $\mathcal P$ the set of all probability measures on $\partial \D.$
Due to the Herglotz representation (\cite{Duren:1983}, Section 1.9) we can write $p(t,z)$ for a. e. $t\geq 0$ as
\begin{equation}
p(t,z) = p_{\mu_t}(z) := \int_{\partial \D} \frac{u+z}{u-z}\, \mu_t(du),
\label{eq:1}
\end{equation}
for some $\mu_t \in \mathcal P.$
 
For $\mu \in \mathcal P$, $\lambda\in\C$ and $w\in\D$ we define the Hamiltonian $H(\mu, \lambda, w)$  by 
$$ H(\mu, \lambda, w) = -\lambda \cdot w \cdot p_\mu(w). $$
Then \eqref{Loewner1} has the form $\dot{w}_t = \frac{\partial}{\partial \lambda} H(\mu_t, \lambda, w(t)).$ \\

Now, if $\{\mu_t\}_{t\geq0}$ leads to an extremal solution $w(t),$ i.e. $w(T) \in \partial R_T(z_0),$ then $\{\mu_t\}_{t\geq0},$ $w(t)$ and $\lambda(t)$ satisfy Pontryagin's maximum  principle; see \cite{MR889459}, p.254, Theorem 3. In our setting we choose complex coordinates and a simple calculation shows that the principle stated in \cite{MR889459} then has the following form:\\
Define $\lambda(t)$ as the solution to the adjoint differential equation 
\begin{equation}\label{costate0}
\dot{\lambda}(t) = -\frac{\partial}{\partial w} H(\mu_t, \lambda(t), w(t)),
\end{equation}
with the initial value condition
\begin{align*}
\lambda(0) = e^{i \beta}, \mbox{ with }\beta\in[0,2\pi).
\end{align*}
Then, for almost every $t\in[0,T]$, we have
\begin{equation}\label{Pontryagin}
\Re H(\mu_t, \lambda(t), w(t)) = \max_{\mu \in \mathcal P} \Re H(\mu, \lambda(t), w(t)),  
\end{equation}
and 
\begin{equation*}
\Re H(\mu_t, \lambda(t), w(t)) = const. \quad \text{for almost all $t\in[0,T].$}
\end{equation*}

In passing we note that equations such as \eqref{Loewner1}, i.e. evolution equations for holomorphic functions, can also be regarded and studied as control systems; see \cite{Roth}.\\  

From \eqref{eq:1} it is easy to see that $\Re H(\mu, \lambda(t), w(t))$ is maximised only for point measures, i.e. when 
$$ H(\mu, \lambda, w) =   -\lambda \cdot w \cdot \frac{u+w}{u-w}$$
for some $u\in \partial\D.$ Thus, for almost every $t\geq 0,$ $H(\mu_t,\lambda(t),w(t))=-\lambda(t)\cdot w(t)\cdot \frac{\kappa(t)+w(t)}{\kappa(t)-w(t)},$ where $\kappa:[0,T]\to \partial \D$ is measurable and \eqref{Loewner1}, \eqref{costate0} become
\begin{equation}
\label{Loewner2} \dot{w}(t) = -w(t) \cdot \frac{\kappa(t)+w(t)}{\kappa(t)-w(t)},\quad w(0)=z_0\in\D, \end{equation}
\begin{equation}\label{costate}
\dot{\lambda}(t) = -\lambda(t)\cdot \frac{w(t)^2-2\kappa(t)w(t)-\kappa(t)^2}{(\kappa(t)-w(t))^2}, \quad \lambda(0)=e^{i\beta}.
\end{equation}
We now optimise the Hamiltonian by rewriting 
\begin{equation*}
\max_{\kappa\in\partial\D}\Re -w \lambda\cdot \frac{\kappa+w}{\kappa-w} =
\max_{\phi\in\R}\Re( -\lambda  w(m+re^{i\phi}))=r|\lambda w|-m\Re(\lambda w),
 \end{equation*}
where
 \begin{align*}
m=\frac{1+|w|^2}{1-|w|^2},\qquad 
r=\frac{2|w|}{1-|w|^2}, \qquad
e^{i\phi}&=\frac{w-|w|^2\kappa}{|w|\kappa-w|w|}.
\end{align*}

The maximum is then obviously taken at
\begin{equation}\label{OptimalDrive}
\phi=\pi-\arg(\lambda w)\quad \Leftrightarrow \quad \kappa=\frac{w}{|w|}\frac{1+|w|e^{i\phi}}{e^{i\phi}+|w|}=w\frac{|\lambda|-\overline{\lambda w}}{|\lambda||w|^2-\overline{\lambda w}}.
 \end{equation}
 Inserting this into the phase equation (\ref{Loewner2}) yields
 \begin{equation*}\label{Loewner3} \dot{w}= -w \left(m+re^{i\phi}\right), 
 \end{equation*}
 or, in polar coordinates,
  \begin{align}
  &\frac{d}{dt}|w|=-|w|(m+r\cos\phi)=-|w|\left(\frac{1+|w|^2-2|w|\cos\left(\arg\lambda+\arg w\right)}{1-|w|^2}\right),\label{LoewnerPolar1}
\\
&\frac{d}{dt}\arg w=-r\sin\phi=-\frac{2|w|\sin\left(\arg\lambda+\arg w\right)}{1-|w|^2},\label{LoewnerPolar2}
 \end{align}
 and the costate equation (\ref{costate}) reads
 \begin{equation*}\label{costate2}
\dot{\lambda}= 
\lambda\left(m+re^{i\phi}+2|w|\frac{|w|+e^{i\phi}(1+|w|^2)+|w|e^{2i\phi}}{(1-|w|^2)^2}\right),
\end{equation*}
which corresponds to
\begin{align*}
\frac{d}{dt}|\lambda|&=|\lambda|\left(m+r\cos\phi+2|w|\frac{|w|+(1+|w|^2)\cos\phi+|w|\cos2\phi
}{(1-|w|^2)^2}\right)=\\
&=|\lambda|\frac{1-|w|^4+2|w|^2-4|w|\cos\left(\arg\lambda+\arg w\right)+2|w|^2\cos\left(2\arg\lambda+2\arg w\right)}{(1-|w|^2)^2},
\end{align*}
\begin{align}
\frac{d}{dt}\arg \lambda&=r\sin\phi+2|w|\frac{|w|\sin(2\phi)+(1+|w|^2)\sin\phi}{(1-|w|^2)^2}= \nonumber\\
&=\frac{4|w|\sin\left(\arg\lambda+\arg w\right)-2|w|^2\sin\left(2\arg\lambda+2\arg w\right)}{(1-|w|^2)^2}.
\label{lambda2}
\end{align}
Now we introduce the variable
  \begin{equation*}
  x:=\cos\left(\arg\lambda+\arg w\right),
  \end{equation*}
  which reduces our system of equations (\ref{LoewnerPolar1}), (\ref{LoewnerPolar2}), (\ref{lambda2}) to 
    \begin{align}\label{pleasebesolvable1}
&\frac{d}{dt}|w|=-|w|\left(\frac{1+|w|^2-2|w|x}{1-|w|^2}\right)
\end{align}
and
\begin{align}\label{pleasebesolvable2}
&\frac{d}{dt}x=-2|w|(1-x^2) \frac{1+|w|^2-2x|w|}{(1-|w|^2)^2}=2\frac{1-x^2}{1-|w|^2}\frac{d|w|}{dt}
  \end{align}
  with the initial value conditions
  \begin{align}\label{init}
  |w(0)|=z_0, \qquad x(0)=x_0:=\cos\beta.
  \end{align}
 For $x_0^2\neq 1$, separation of variables solves (\ref{pleasebesolvable2}), (\ref{init}) as
      \begin{align*}
      x(t)&=\Phi^{-1}\left(2\arctanh |w(t)|-2\arctanh z_0\right),
                  \end{align*}
 where 
 \begin{align*} 
      \Phi(y)&:=\arctanh y-\arctanh x_0 ,
\end{align*}
which means 
      \begin{align}\label{xsolved}
x(t)&=\tanh\left(2\arctanh |w(t)|+\arctanh x_0-2\arctanh z_0\right)=\nonumber \\ 
&=\frac{(1+|w(t)|^2)\left(x_0-2z_0+x_0z_0^2\right)+2|w(t)|\left(1-2x_0z_0+z_0^2\right)}
{(1+|w(t)|^2)\left(1-2x_0z_0+z_0^2\right)+2|w(t)|\left(x_0-2z_0+x_0z_0^2\right)}=\nonumber \\
&=\frac{(1+|w(t)|^2)A+2|w(t)|B}
{(1+|w(t)|^2)B+2|w(t)|A}
\end{align}
with
\begin{align*}
A&:=x_0-2z_0+x_0z_0^2, \\
B&:=1-2x_0z_0+z_0^2.
\end{align*}
Note that, in fact, the denominator in (\ref{xsolved}) never equals zero for any $x_0\in[-1,1]$, since we have
\begin{align*}
(1+|w|^2)B+2|w|A=0 \quad \Leftrightarrow\quad |w|&=-\frac{A}{B}\pm\frac{\sqrt{A^2-B^2}}{B}\\
&=-\frac{A}{B}\pm\frac{\sqrt{(x_0^2-1)(1-z_0^2)^2)}}{B},
\end{align*}
which only yields real terms for $x_0^2=1$, and in this case the only solution is
\begin{align*}
|w|=-\frac{A}{B}=\pm1\not\in(0,1).
\end{align*}
Therefore, (\ref{xsolved}) is for all $x_0\in[-1,1]$ the solution to the initial value problem (\ref{pleasebesolvable2}), and thus (\ref{pleasebesolvable1}) can be simplified to
      \begin{align*}
      \frac{d}{dt}|w(t)|=-|w|\frac{B(1-|w|^2)}{B(1+|w|^2)+2|w| A}, \quad |w(0)|=z_0.
      \end{align*}
      The function
       \begin{align*}
      \Psi(y)&:=(A + B) \log(1 - y) - B \log(y) - (A - B) \log(1 + y)
      \end{align*}
      is strictly monotonous on the interval $(0,1)$, since its derivative is zero-free. Hence it is invertible, and 
      \begin{align*}
|w(t)|=\Psi^{-1}(B t +\Psi(z_0)),
\end{align*}
is the solution to the initial value problem (\ref{pleasebesolvable1}), which can be verified by calculation. \\
To determine the value set $R_T(z_0)$,
we solve the remaining initial value problem (\ref{LoewnerPolar2}), which now reads
\begin{align*}
\frac{d}{dt}\arg w(t)=\pm\frac{2\sqrt{B^2-A^2}}{B(1+|w(t)|^2)+2A|w(t)|}, \quad \arg w(0)=0.
\end{align*}
If we write
\begin{align*}
\arg w(t)=-G(|w(t)|),
\end{align*}
where $G$ is the solution to
\begin{align*}
\frac{d}{d|w|}G(|w|)=\frac{2\sqrt{B^2-A^2}}{B(1-|w|^2)}, \quad G(0)=0,
 \end{align*}
 then
\begin{align*}
\arg w(t)&=\frac{\pm 2\sqrt{B^2-A^2}}{B}\left(\arctanh z_0-\arctanh|w(t)|\right).
\end{align*} 

We can therefore describe candidates for the boundary points of the set $R_T(z_0)$ as follows:\\
For $x_0\in[-1,1]$, let $r=r(T,x_0)$
be the (unique) solution to the equation
\begin{align}\label{kai}
(1+x_0)(1-z_0)^2\log(1-r)+(1-x_0)(1+z_0)^2\log(1+r)-(1-2x_0z_0+z_0^2)\log r = \nonumber\\
(1+x_0)(1-z_0)^2\log(1-z_0)+(1-x_0)(1+z_0)^2\log(1+z_0)-(1-2x_0z_0+z_0^2)\log e^{-T
}z_0,
\end{align}
then $\partial R_T(z_0)$ consists of a subset of the two curves

\begin{align*}
C_{\pm}(z_0)=\left\{w_{\pm}(x_0)= r(T, x_0)e^{\pm i\sigma(T, x_0)}\,:\,  x_0\in[-1,1]\right\},
\end{align*}
where $$\sigma(T, x_0)=\frac{2(1-z_0^2)\sqrt{1-x_0^2}}{1-2x_0z_0+z_0^2}\left(\arctanh z_0 - \arctanh r(T, x_0)\right). $$

First we consider the function $x_0 \mapsto r(T,x_0)$: By solving \eqref{kai} for $T$ and then taking the derivative  with respect to $x_0$, we obtain

$$\frac{\partial}{\partial x_0}r(T,x_0)=
- \frac{(1-z_0)^2 r(T,x_0)(1-r^2(T,x_0))\left(\log\left(\frac{1+r(T,x_0)}{1-r(T,x_0)}\right)-\log\left(\frac{1+z_0}{1-z_0}\right)\right)}{B(B(1+r^2(T,x_0))+2A\phantom{\cdot}r(T,x_0))},$$
and since the only zeros of this term lie at $r(T,x_0)=0$, $r(T,x_0)=\pm1$ and $r(T,x_0)=z_0$, this immediately shows that $x_0\mapsto r(T,x_0)$ is strictly increasing.\\
In particular, the curves $C_+(z_0)$ and $C_-(z_0)$ do not hit themselves.\\

Now we consider the first case where $z_0<\tanh\frac{\pi}{2}$. Here, the curves never hit the negative real axis:\\
As the function
$$x_0\mapsto \frac{2(1-z_0^2)\sqrt{1-x_0^2}}{1-2x_0z_0+z_0^2}$$
reaches its single maximal value 2 at $x_0=\frac{2z_0}{1+z_0^2}$, we have
$$\sigma(T,x_0)=\frac{2(1-z_0^2)\sqrt{1-x_0^2}}{1+2x_0z_0+z_0^2}\left(\arctanh z_0-\arctanh r(T,x_0)\right)< 2\cdot (\pi/2-0) = \pi.$$
Thus, they intersect only on the positive real axis and, as $\sigma(T,x_0)=0$ if and only if $x_0=\pm1$, this happens exactly at $x_0=\pm1$. Hence, the full set $C_+(z_0)\cup C_-(z_0)$ forms the boundary of $R_T(z_0)$. Since $R_T(z_0)$ is obviously bounded, it has to consist of the bounded region enclosed by the two curves.\\

Next assume that $z_0>\tanh\frac{\pi}{2}$. We have
$$\frac{\partial}{\partial x_0}\sigma(T,x_0)=
-\frac{1-z_0^2}{B\sqrt{1-x_0^2}}\left(\arctanh z_0-\arctanh r(T,x_0)\right)\frac{A(1+r(T,x_0)^2)+2B r(T,x_0)}{B(1+r(T,x_0)^2)+2A r(T,x_0)}.$$
The zeros of this term lie clearly at the points $x_0\neq\frac{2z_0}{1+z_0^2}$ with
$$r(T,x_0)=\frac{-B\pm\sqrt{B^2-A^2}}{A}.$$
Since $$
\frac{-B-\sqrt{B^2-A^2}}{A}\qquad \left\{ 
\begin{array}{ll}
\geq 1 & \mbox{ for }x_0<\frac{2z_0}{1+z_0^2}\\
<0& \mbox{ for }x_0>\frac{2z_0}{1+z_0^2},
\end{array}
\right.$$
it is clear that this term can be ignored. We focus on the equality
\begin{align}\label{rboese}   
r(T,x_0)=\frac{-B+\sqrt{B^2-A^2}}{A}
\end{align}
and note that here the term on the right-hand side is well-defined for all $x_0\in[-1,1]$, and strictly decreasing on this interval, taking values between $-1$ and $1$. Therefore, $x_0\mapsto h(x_0):=\frac{-B+\sqrt{B^2-A^2}}{A} - r(T,x_0)$  is continuous on $[-1,1]$, strictly decreasing, and we have $h(-1)\geq0$ and $h(1)\leq0$.
Thus (\ref{rboese}) has exactly one solution $x_0=x^*$ on $[-1,1]$, and the function $x_0\mapsto\sigma(T,x_0)$ increases from 0 to $\sigma(T,x^*)$ and decreases again to 0.\\

If $T$ is so small that equation \eqref{split} has no solution, then we are again in the same situation: the two curves intersect only twice, namely for $x_0=\pm1$, and $R_T(z_0)$ is the closed region bounded by the two curves.\\

\begin{figure}[h]
\begin{minipage}{0.5\textwidth}
\includegraphics[width=0.8\textwidth]{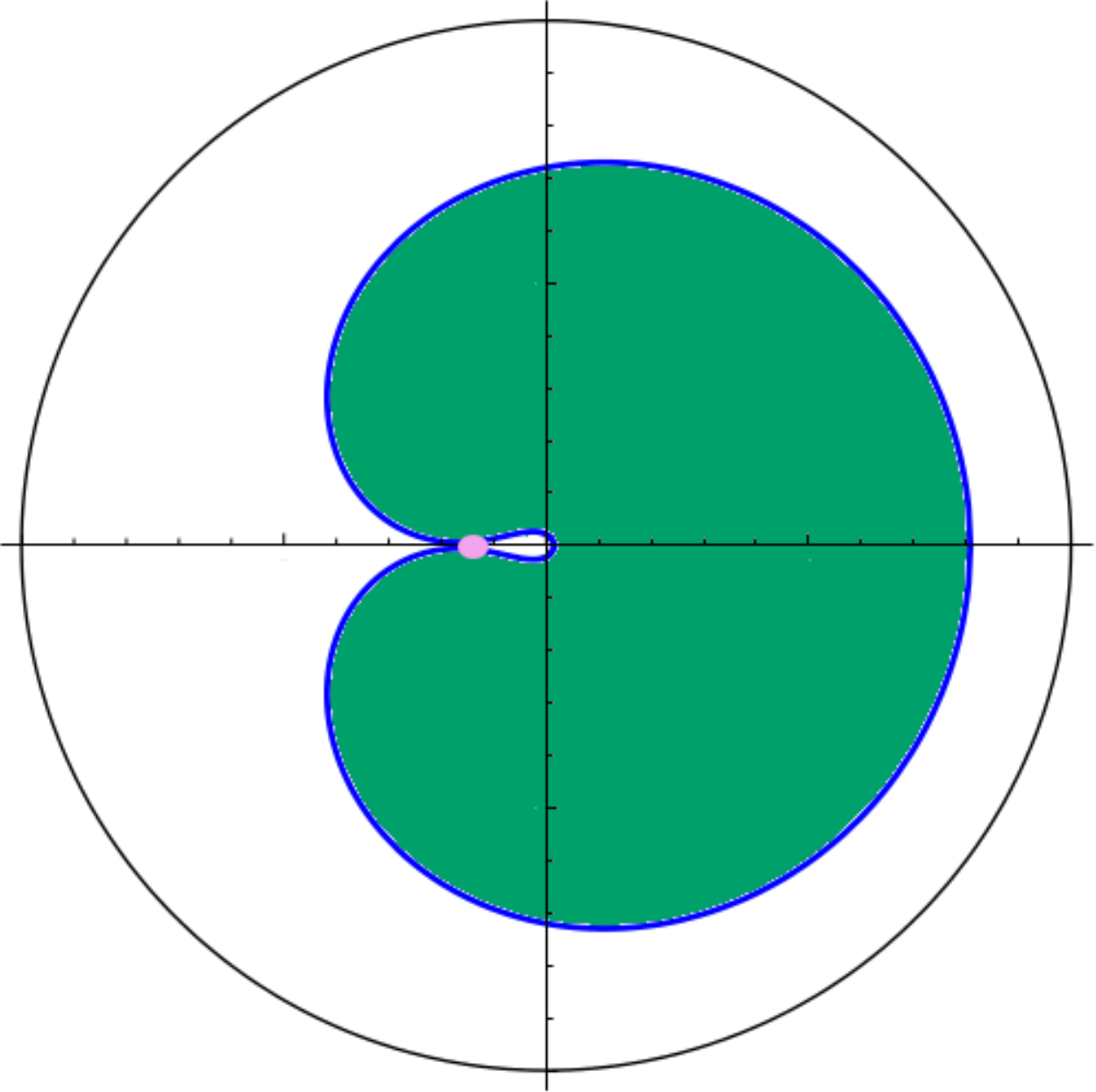}
\caption{$V_{T^*}(z_0)$}
\label{firstcontact}
\end{minipage}
\hfill
\begin{minipage}{0.5\textwidth}
\includegraphics[width=0.8\textwidth]{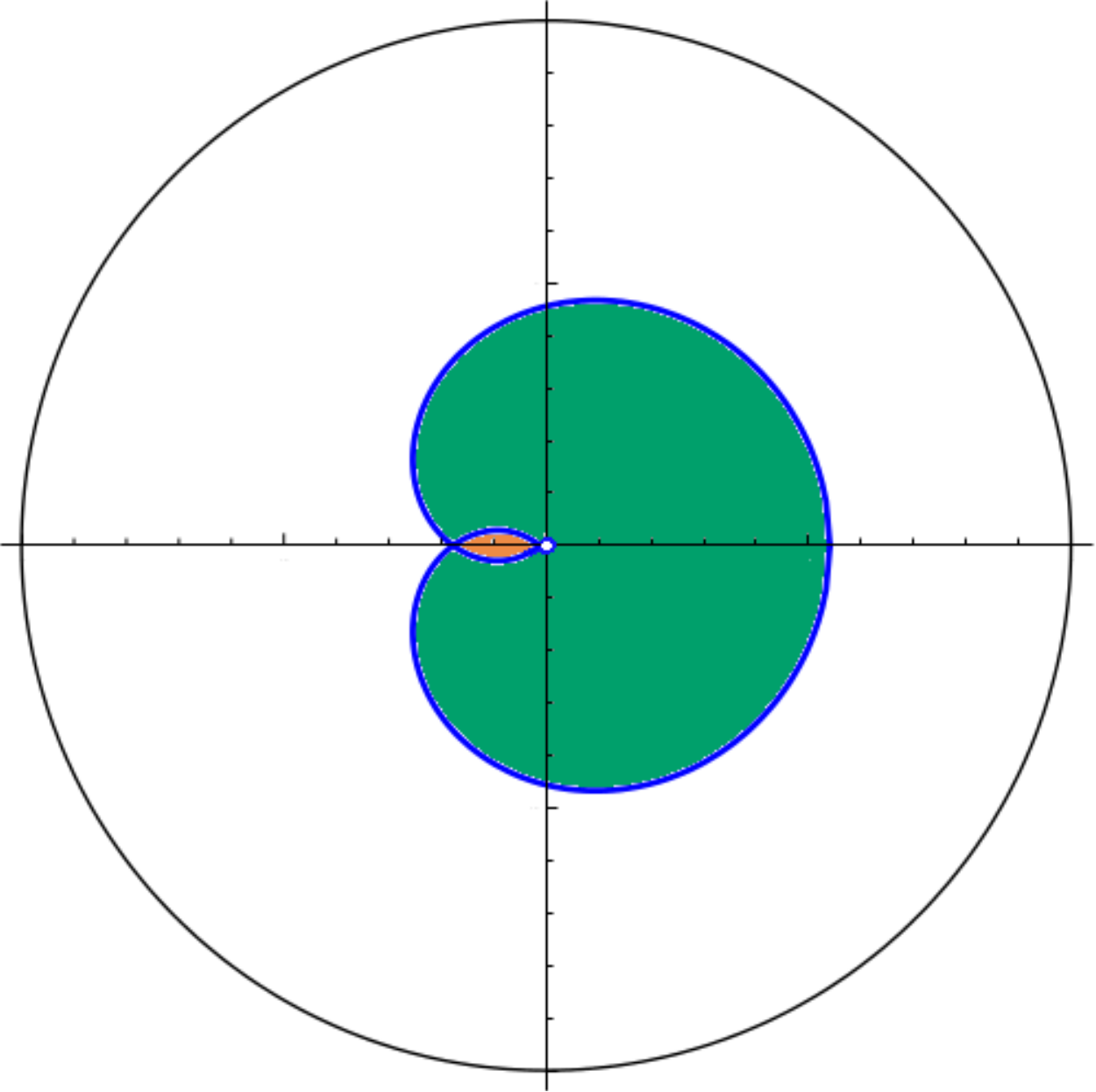}
\caption{$V_{T^*+\varepsilon}(z_0)$}
\label{decompV}
\end{minipage}
\begin{center}
The evolution of the decomposition of $\D$ by $C_\pm(z_0)$
\end{center}
\end{figure}

There is a $T^*$ such that \eqref{split} admits a solution, but has no solution for any $T<T^*$. At this $T^*$, the curves $C_\pm(z_0)$ will meet for the first time, i.e. $\sigma(T^*,x^*)=\pi$. This means that at $x^*$, the curves both touch $\R^-$ at some point $z^*$, see Figure \ref{firstcontact}, and $R_T(z_0)$ (shown in green) is no longer simply connected, since the component containing the origin can obviously not be part of $R_T(z_0)$. \\
 
For slightly larger $T$, the curves $C_\pm(z_0)$ intersect on $\R^-$ twice and $\D\setminus (C_+(z_0)\cup C_-(z_0))$ has four components, see Figure \ref{decompV}. We denote by $K_T(z_0)$ the component (shown in orange) that arises from the intersection of the two curves near $x_0=x^*.$
Obviously, the component that contains the origin, as well as the "exterior" component (both shown in white) cannot be part of $R_T(z_0)$. For reasons of continuity, the "large interior" component (shown in green) must belong to $R_T(z_0)$. 
It remains to show that $K_T(z_0)$ also belongs to $R_T(z_0)$:\\
Since $z^*=w(T^*)$ for a solution $w(t)$ of the Loewner equation (\ref{Loewner2}), we know that $R_T(z_0)$ contains the set $R_{T-T*}(z^*)$, which we determined already if $T-T^*$ is small enough. In particular, $R_T(z_0)$ contains infinitely many points of $\R^-.$ If $K_T(z_0)$ was not included in $R_T(z_0)$, then $R_T(z_0)\cap \R^-$ would consist of only two points, a contradiction. \\
For reasons of continuity, the set $R_T(z_0)$ will have the form described in the theorem for any larger $T$ as well, and this concludes the proof. 

\begin{remark}
If $w_0 \in \partial V_T(z_0),$ then there exists exactly one control function $\kappa(t)$ such that the solution $\{f_t\}_{t\in[0,T]}$ of \eqref{Loewner0} with  $p(t,z)=\frac{\kappa(t)+z}{\kappa(t)-z}$ satisfies $f_T(z_0)=w_0.$ Equation \eqref{OptimalDrive} shows that $\kappa(t)=\exp(i\varphi(t))$  is continuously differentiable. From \cite{MarshallRohde:2005}, Theorem 1.1, it follows that $f$ is a slit mapping in this case, i.e. $f$ maps $\D$ conformally onto $\D\setminus \gamma$, where $\gamma$ is a simple curve.
\end{remark}

\section{Value sets for the inverse functions}\label{sec:inverse}

Firstly, in analogy to \cite{MR3262210} and the set $\mathcal V(z_0)$, we describe the set $$\mathcal W(z_0):=\{f^{-1}(z_0)\, :\, f\in \mathcal{S} \text{ with $z_0\in f(\D)$}\}.$$

In the following we write $d_\D(0,z),$ $z\in\D$, for the hyperbolic distance between $0$ and $z$ (using the hyperbolic metric with curvature -1), i.e. 
$d_\D(0,z) = 2\arctanh(|z|) = \log\left(\frac{1+|z|}{1-|z|}\right).$

\begin{theorem}\label{th3}
We have
\begin{align*}
\mathcal W(z_0)&=\{f^{-1}(z_0)\, : \, f:\D\to\D \text{ univalent, } f(0)=0,\phantom{i} f'(0)>0 \text{ with $z_0\in f(\D)$} \}\\
&=
\{re^{i\sigma}\,:\,d_\D(0,r)\geq |\sigma|+d_\D(0,z_0),\, \sigma \in [-\pi,\pi]\}.
\end{align*}
\end{theorem}

\begin{figure}[h]
\begin{center}
\includegraphics[width=6cm]{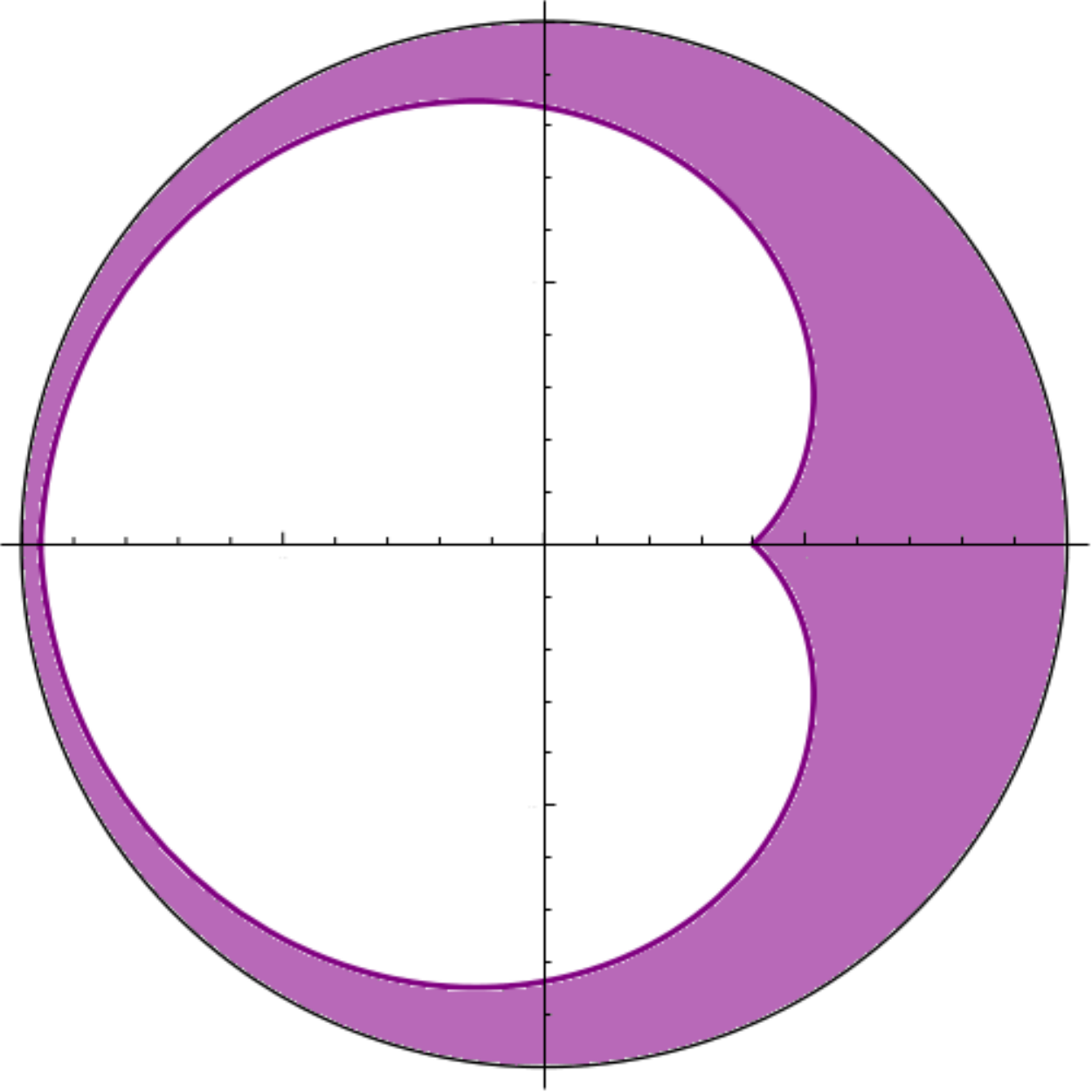}
\end{center}
\caption{The set  $\mathcal W(0.4)$ }
\end{figure}

Furthermore, we will determine the value set 
$$W_T(z_0):=\{f^{-1}(z_0)\, :\, f\in \mathcal{S}_T \text{ with $z_0\in f(\D)$}\}$$ for the inverse functions:
\begin{theorem}\label{th2} Let $z_0\in(0,1).$
For $x_0\in[-1,1)$ and $T>0$, let $r=r(T,x_0)$
be the (unique) positive solution to the equation
\begin{align*}
(1-x_0)(1-z_0)^2\log(1-r)+(1+x_0)(1+z_0)^2\log(1+r)-(1+2x_0z_0+z_0^2)\log r=\\
(1-x_0)(1-z_0)^2\log(1-z_0)+(1+x_0)(1+z_0)^2\log(1+z_0)-(1+2x_0z_0+z_0^2)\log e^{T}z_0
\end{align*}
and let 
$$\sigma(T,x_0)=\frac{2(1-z_0^2)\sqrt{1-x_0^2}}{1+2x_0z_0+z_0^2}\left(\arctanh r(T,x_0)-\arctanh z_0\right).$$
If 
$$T< T^*:=\log\frac{(1+z_0)^2}{4z_0},$$
then $r(T,x_0)$ can be extended continuously to $x_0=1$ and we have $W_T(z_0)=\overline{W_T(z_0)}\subset \D$, and $W_T(z_0)$ is the closed region bounded by the two curves
\begin{align*}
D_\pm(z_0)&:=\left\{ r(T,x_0)e^{\pm i\sigma(T,x_0)}\,:\,  x_0\in[-1,1]\right\}.
\end{align*}
Now let $T \geq T^*$ and define the two curves
\begin{align*}
\widetilde D_\pm(z_0)&:=\left\{ r(T,x_0)e^{\pm i\sigma(T,x_0)}\,:\,  x_0\in[-1,1)\right\}.
\end{align*} Here we have two cases: if $T$ is small enough that $\widetilde D_+(z_0)$ and $\widetilde D_-(z_0)$ intersect only at $x_0=-1$, then
$\overline{W_T(z_0)}$ intersects $\partial\D$ and $\overline{W_T(z_0)}$ is bounded by the two curves
$\widetilde D_\pm(z_0)$
and by the part of $\partial\D$ between the intersection points with the curves which includes the point $1$. \\

Otherwise, the two curves intersect on $\R^-$ for the first time for some $x_0=\chi\in(-1,1)$ and $\overline{W_T(z_0)}$ is the closed region bounded by $\partial \D$ and the two curves
$$
\widehat D_\pm(z_0):=\left\{ r(T,x_0)e^{\pm i\sigma(T,x_0)}\,:\,  x_0\in[-1,\chi]\right\}.
$$
In the last two cases we obtain $W_T(z_0)$ from $W_T(z_0)=\overline{W_T(z_0)}\cap \D$.
\end{theorem}

\begin{figure}[h]
\begin{center}
\includegraphics[width=6cm]{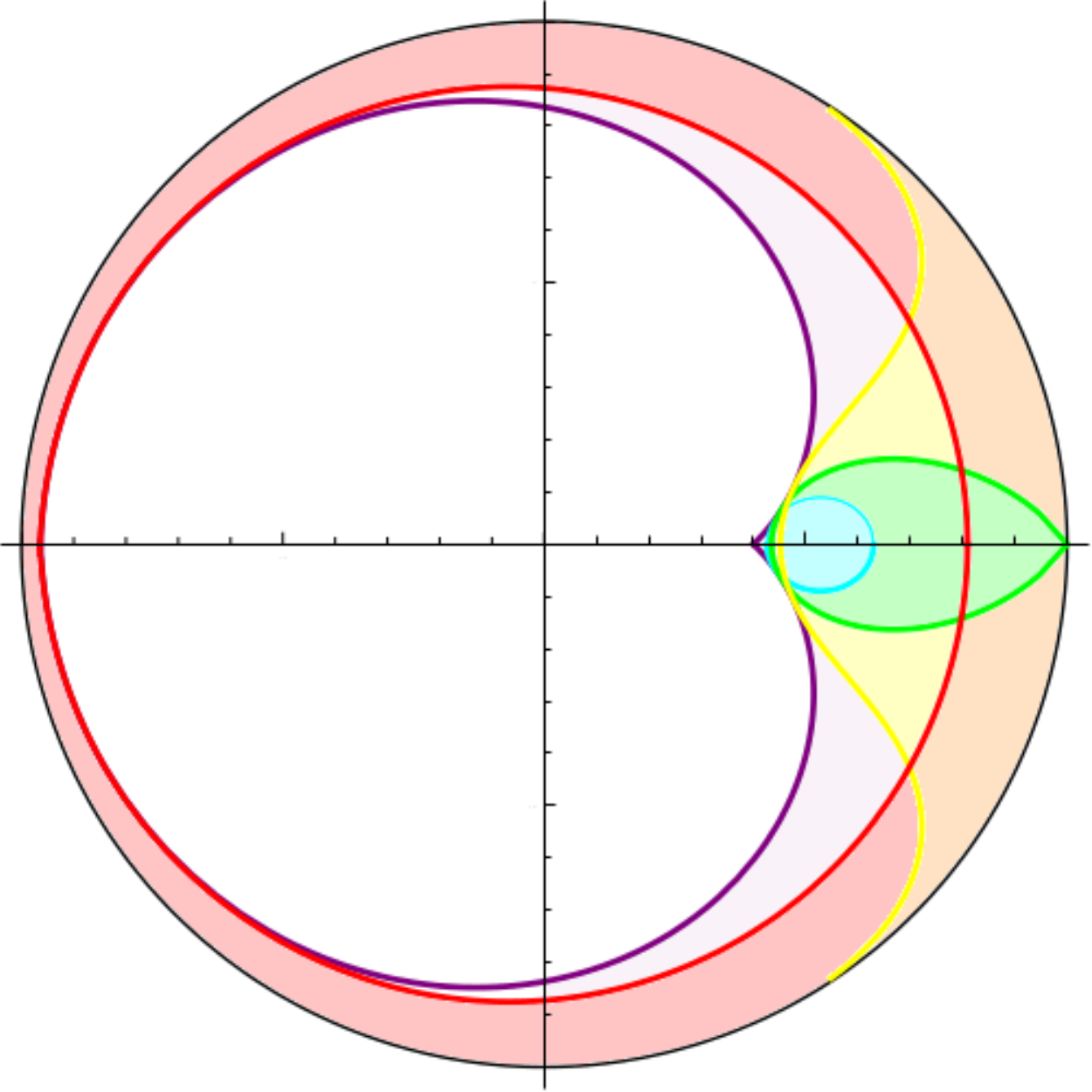}
\end{center}
\caption{$W_T(0.4)$ for $T=0.15, T^*,0.3, 3.$  }
\end{figure}

\section{Proofs of Theorem \ref{th3} and  \ref{th2} }\label{sec:proof2}
The proof of Theorem \ref{th2} is analogous to that of Theorem \ref{thm:1}: we consider the inverse Loewner equation
\begin{align}\label{Loewnerinvers0}
\dot{w}(t) = w(t) \cdot p(t,w(t)),\quad w(0)=z_0\in\D,
\end{align}
where $p(t,z)$ is a Herglotz function.

Here, a solution $t\mapsto w(t)$ may not exist for all time, i.e. there might be a $t_{max}>0$ such that $w(t)\in\D$ for all $t<t_{max}$ but $|w(t)|\to1$ for $t\uparrow t_{max}.$ In this case, the (classical) solution to \eqref{Loewnerinvers0} ceases to exist at $t_{max}$. We define the reachable set 
$$R'_T(z_0)=\{w(T)\,:\, w:[0,T]\to\D \text{ solves \eqref{Loewnerinvers0}}\}.$$
Note that we assume here that $w(t)$ exists up to $t=T$ and $w(T)\in\D.$ \\
Then $W_T(z_0)=R'_T(z_0)$ is closed in the relative topology on $\D$, and we have
$$ W_T(z_0)  = \overline{W_T(z_0)} \cap \D.$$

Next we describe the boundary $\partial R'_T(z_0)$ by applying the maximum principle to \eqref{Loewnerinvers0}. For $\mu \in \mathcal P$, $\lambda\in\C$ and $w\in\D$ we now have the Hamiltonian 
$$ H'(\mu, \lambda, w) = \lambda \cdot w \cdot p_\mu(w). $$
Since the only difference to the case $R_T(z_0)$ consists in the sign of the left hand side of the Loewner differential equation, we can use the exact same ideas as above. Equation \eqref{Loewnerinvers0} reduces to 
\begin{align}\label{Loewnerinvers}
\dot{w}(t) = w(t) \cdot \frac{\kappa(t)+w(t)}{\kappa(t)-w(t)},\quad w(0)=z_0\in\D,
\end{align}
where $\kappa:[0,T]\to\partial \D$ is measurable. The condition \eqref{OptimalDrive} that is satisfied by trajectories leading to boundary points now corresponds to
$$\phi=-\arg(\lambda w),$$
which means we have to solve the system of equations
\begin{align}
\label{s:l}\frac{d}{dt}|w|=|w|\frac{1+|w|^2+2|w| x}{1-|w|^2},\quad |w(0)|=z_0,\\
\nonumber\frac{d}{dt}x=-2\frac{1-x^2}{1-|w|^2}\frac{d|w|}{dt}, \quad x(0)=:x_0\in [-1,1].
\end{align}
We are left with
\begin{align*}
 x(t)&=\Delta^{-1}\left(2\arctanh |w(t)|-2\arctanh z_0\right),
\end{align*}
where
$$\Delta(y)=\arctanh x_0-\arctanh y,$$
and thus 
\begin{align*}
x(t)&=\tanh\left(\arctanh x_0+2\arctanh z_0-2\arctanh |w(t)|\right)=\\
&=\frac{(1+|w|^2)G-2H |w|}{(1+|w|^2)H-2G|w|},
\end{align*}
where
\begin{align*}
G&:=x_0+2z_0+x_0z_0^2,\\
H&:=1+2x_0z_0+z_0^2.
\end{align*}
Note that, again, this last term for $x$ is valid for any $x_0\in[-1,1]$. \\
We hence arrive at
\begin{align*}
\frac{d|w|}{dt}=\frac{H|w|(1-|w|^2)}{H(1+|w|^2) - 2 G|w|},
\end{align*}
or
\begin{align*}
|w(t)|=\Theta^{-1}\left(-H t +\Theta(z_0)\right)
\end{align*}
with
$$\Theta(y)=(H-G)  \log(1-y)-H \log y+(G+H)\log(1+y).$$
The differential equation for the argument of the optimal trajectory $w$ reads
\begin{align*}
\frac{d}{dt}\arg w(t)=\pm\frac{2|w|\sqrt{H^2-G^2}}{(1+|w|^2)H-2G|w|},
\end{align*}
which means
\begin{align*}
\arg w(t)=\pm\frac{2\sqrt{H^2-G^2}}{H}\left(\arctanh|w|-\arctanh z_0\right).
\end{align*}

 We can now describe the sets $R'_T(z_0)$:\\
Let $x_0\in[-1,1).$ Then $\Theta\left((0,1)\right)=(-\infty,\infty)$ and $\Theta$ is strictly decreasing. Thus there is exactly one solution $r=r(T,x_0)$
of the equation
\begin{align}\label{kaihe}
(1-x_0)(1-z_0)^2\log(1-r)+(1+x_0)(1+z_0)^2\log(1+r)-(1+2x_0z_0+z_0^2)\log r= \nonumber \\
(1-x_0)(1-z_0)^2\log(1-z_0)+(1+x_0)(1+z_0)^2\log(1+z_0)-(1+2x_0z_0+z_0^2)\log e^{T
}z_0. \end{align}
Furthermore we define the two curves
\begin{align*}
\widetilde D_\pm(z_0)&:=\left\{  r(T,x_0)e^{\pm i\sigma(T,x_0)}\,:\,  x_0\in[-1,1)\right\},
\end{align*}
where 
$$\sigma(T,x_0)=\frac{2(1-z_0^2)\sqrt{1-x_0^2}}{1+2x_0z_0+z_0^2}\left(\arctanh r(T,x_0)-\arctanh z_0\right).$$

We take a closer look at the absolute value $r(T,x_0)$. \\
Firstly, the function $x_0\mapsto r(T,x_0)$ is strictly increasing:\\
By solving (\ref{kaihe}) for $T$ and then deriving with respect to $x_0$, we can calculate
$$\frac{\partial}{\partial x_0}r(T,x_0)=  \frac{(1-z_0)^2 r(T,x_0)(1-r^2(T,x_0))\left(\log\left(\frac{1+r(T,x_0)}{1+z_0}\right)-\log\left(\frac{1-r(T,x_0)}{1-z_0}\right)\right)}{H(H(1+r^2(T,x_0))-2G\phantom{\cdot}r(T,x_0))},$$
and since the only zeros of this term lie at $r(T,x_0)=0$, $r(T,x_0)=\pm1$ and $r(T,x_0)=z_0$, this immediately shows that $x_0\mapsto r(T,x_0)$ is strictly increasing in $x_0$ for $T>0$. \\
Hence, we can define $r(T,x_0)$ also for $x_0=1$.\\

Note that for $x_0=1$, (\ref{kaihe}) simplifies to
$$2\log(1+r)-\log r=2\log(1+z_0)-\log z_0 -T,$$
which means that the curves $D_+(z_0)$ and $D_-(z_0)$ will hit the boundary of the unit circle for the first time for
$$T=T^*:=\log\frac{(1+z_0)^2}{4z_0}.$$

Next we take a closer look at the behaviour of the argument $\sigma(T,x_0)$ of the curve. We calculate 
\begin{align*}
\frac{\partial}{\partial x_0}\sigma(T,x_0)=&
\frac{2(1-z_0^2)\left(\arctanh r(T,x_0)-\arctanh z_0\right)}{H^2}
\left( \frac{2r(T,x_0)\sqrt{1-x_0^2}(1-z_0^2)^2}{(H(1+r^2(T,x_0))-2G\phantom{\cdot}r(T,x_0))}
-
\frac{G}{\sqrt{1-x_0^2}}
\right).
\end{align*}
Since
\begin{align*}
H(1+r^2(T,x_0))-2G\phantom{\cdot}r(T,x_0)\geq0 \text{ for all $x_0\in(-1,1)$, $z_0\in(0,1)$ and $r(T,x_0)\geq z_0$},
\end{align*}
the term is non-negative if and only if
\begin{align*}
2r(T,x_0)(1-x_0^2)(1-z_0^2)^2\geq (HG(1+r^2(T,x_0))-2G^2\phantom{\cdot}r(T,x_0)),
\end{align*}
or
\begin{align*}
H(G-2H\cdot r(T,x_0)+G\cdot r^2(T,x_0))\leq 0,
\end{align*}
which is equivalent to
\begin{align}\label{sigma'}
\frac{H-\sqrt{H^2-G^2}}{G}\leq r(T,x_0)\leq \frac{H+\sqrt{H^2-G^2}}{G}
\end{align}
The inequality to the right always holds, since 
$$\frac{H+\sqrt{H^2-G^2}}{G} \quad\left\{ 
\begin{array}{ll}
\leq0 &\mbox{ for }x_0<-\frac{2z_0}{1+z_0^2},\\
>1 &\mbox{ for }x_0>-\frac{2z_0}{1+z_0^2},\\
\end{array} \right.$$
and of course
$$0<r(T,x_0)\leq1 \mbox{ for all } x_0\in[-1,1).$$

The curves $\widetilde D_+(z_0)$ and $\widetilde D_-(z_0)$ can only intersect on $\R$, i.e. $\sigma(T,x_0)=k\cdot \pi.$ 
Obviously, $\sigma(T,x_0)\geq0$ for all $x_0$ so that $k\geq0$ when the two curves intersect.\\

Next we show that \begin{equation}\label{grow}
\text{$\frac{\partial}{\partial x_0}\sigma(T,x_0)>0$\quad if \quad $\sigma(T,x_0)\geq \pi.$}
\end{equation} 
We have
\begin{align*}
\log\left(1+\frac{2H-2\sqrt{H^2-G^2}}{G-H+\sqrt{H^2-G^2}}\right) \leq \frac{2H-2\sqrt{H^2-G^2}}{G-H+\sqrt{H^2-G^2}}
\leq \frac{\pi H}{\sqrt{H^2-G^2}},
\end{align*}
for $$2(H\sqrt{H^2-G^2}-H^2+G^2)\leq \pi (H\sqrt{H^2-G^2}- H^2+HG) ,$$
and thus
$$
r(T,x_0)>\tanh\left(\frac{\pi H}{2(1-z_0^2)\sqrt{1-x_0^2}}\right)  \geq \frac{H-\sqrt{H^2-G^2}}{G} .$$
Thus, \eqref{sigma'} is satisfied in this case and $\frac{\partial}{\partial x_0}\sigma(T,x_0)>0$.\\

Now we consider the first case $T<T^*:$\\
Here, $r(T,1)<1$ and $\sigma(T,x_0)$ is defined also for $x_0=1.$ Furthermore, $\sigma(T,\pm1)=0,$ i.e. the two curves  $D_+(z_0)$ and $D_-(z_0)$ intersect for $x_0=\pm1$ on the positive real axis. Assume that the curves intersect more than twice. As $\sigma(T,x_0)>0$ for all $x_0\in(-1,1)$ there must be some $\rho\in(-1,1)$ with $\sigma(T,\rho)=\pi.$ This is a contradiction: the function $x_0\mapsto \sigma(T,x_0)$ is increasing for $x_0\in[\rho,1]$ because of \eqref{grow}, but $\sigma(T,1)=0.$
Thus, the two curves don't intersect for $x_0\in(-1,1).$\\ 
Consequently, the set $R'_T(z_0)$ is the closed region enclosed by 
$D_+(z_0) \cup D_-(z_0)$.\\

Next let $T=T^*.$ Then $\overline{R'_{T^*}(z_0)}$ is still the closed region bounded by 
$D_+(z_0) \cup D_-(z_0)$, but $R'_{T^*}(z_0)=\overline{R'_{T^*}(z_0)}\setminus\{1\}$ is not closed anymore.\\
In passing we note that it is not difficult to show that the solution $w(t)$ of \eqref{Loewnerinvers} with $\kappa(t)\equiv 1$ satisfies $\lim_{t\to T^*}w(t)=1$ and that this case corresponds to a mapping $f\in \mathcal S_{T^*}$ that maps $\D$ onto $\D$ minus the slit $[z_0,1].$]\\

Now let $T>T^*$.\\ 
It is easy to see that the function $\Theta$, which defines $r(T,x_0)$, is strictly decreasing, and that therefore, for fixed $x_0$, the term $r(T,x_0)$ is strictly increasing with growing $T$. Thus we know that we still have
%Furthermore, $\Theta\left((0,1)\right)=(-\infty,\infty)$ for $x_0\neq1$, and  $\Theta\left((0,1)\right)=(2\log2,\infty)$ for $x_0=1$.\\
$$r(T,x_0)\to1 \text{ for }x_0\to1.$$
The driving function $\kappa(t)\equiv1$ will now generate a mapping from $\D$ onto $\D\setminus [a,1]$ with $a<z_0.$ From this it is easy to deduce that $$L(T):=\liminf_{x_0 \to 1} \sigma(T, x_0)>0.$$
Furthermore, $L(T)$ is increasing in $T\in[T^*,\infty)$: For a point $p=e^{i\alpha}\in\partial \D$ the driving function $\kappa(t)\equiv -e^{i\alpha}$ has the property that $-p\cdot \frac{\kappa(t)+p}{\kappa(t)-p} = 0.$ Thus, if $e^{i\alpha}\in \overline{R'_T(z_0)},$ then also $e^{i\alpha}\in \overline{R'_S(z_0)}$ for all $S\geq T.$\\

If $T$ is so small that $L(T)\leq \pi,$ then the curves $\widetilde D_\pm(z_0)$ do not intersect in $\D$ a second time besides $x_0=-1$ for the same reason as in the case $T<T^*.$ Here, $\overline{R'_T(z_0)}$ is the closed region which is bounded by $\widetilde D_+(z_0)$ and $\widetilde D_-(z_0)$ and the part of $\partial\D$ which includes the point $1$.\\

\begin{figure}[h]
\begin{center}
\includegraphics[width=6cm]{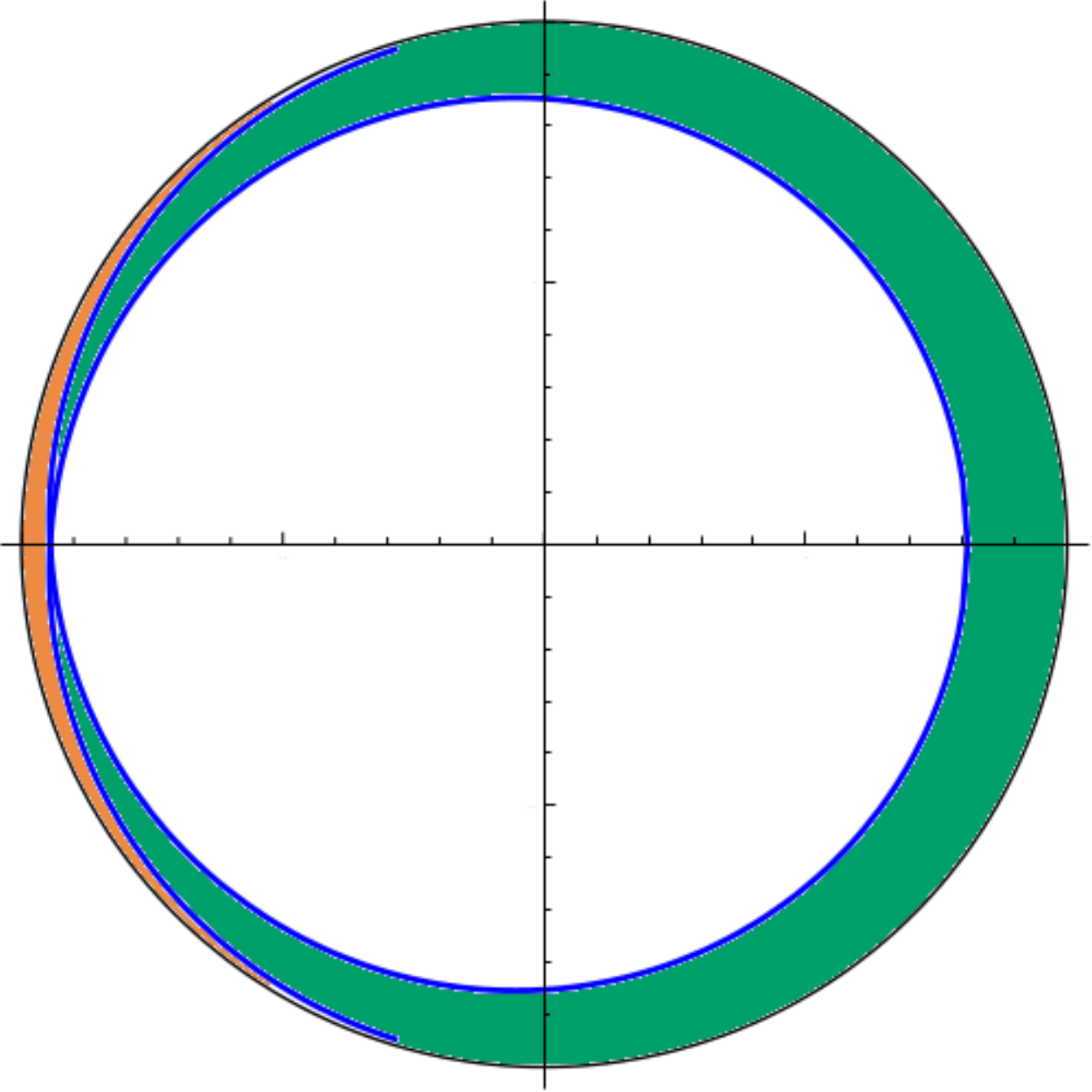}
\end{center}
\caption{The decomposition of $\D$ by  $\widetilde D_\pm(z_0)$}
\label{decomp}
\end{figure}

Finally, let $L(T)>\pi.$ The curves
$\widetilde D_\pm(z_0)$ will meet at $x_0=-1$, and then intersect again on the negative real axis before hitting $\partial \D$. Because of \eqref{grow} they don't intersect more than twice provided that $T>T^*$ is small enough.
Hence, in this case, $\D\setminus (\widetilde D_+(z_0) \cup \widetilde D_-(z_0))$ has three components, see Figure \ref{decomp}.

There is a simply connected component that is bounded by $\widehat D_+(z_0)\cup\widehat D_-(z_0)$ and does not touch $\partial\D$, and  two simply connected components that do touch the boundary $\partial\D$. We denote by 
$W^\pm_T(z_0)$ the components that touch the points $+1$ (shown in green), or, respectively, $-1$ (shown in orange). It is clear that $\overline{R'_T(z_0)}$ has to consist of either $W^+_T(z_0)$, or $W^-_T(z_0)$, or the union of both. If it were equal to only one of the sets $W^\pm_T(z_0)$, this would imply that $\overline{R'_T(z_0)}$ is bounded away from parts of $\partial\D$, although $\overline{R'_t(z_0)}$, with some $t<T$, already touched these segments of $\partial\D$ - a contradiction. Thus, we must have 
$\overline{R'_T(z_0)}=\overline{W^+_T(z_0) \cup W^-_T(z_0)}$, and thus $\overline{R'_T(z_0)}$ is exactly the closed region bounded by $\partial\D$ and (in the interior) by $\widehat D_+(z_0)\cup\widehat D_-(z_0)$.\\
The same consideration applies as well for the case of more than {two intersections} of $\widetilde D_\pm(z_0)$ with $\R$, and for reasons of continuity, the inner boundary of $\overline{R'_T(z_0)}$ has to consists of $\widehat D_+(z_0)\cup\widehat D_-(z_0)$ in these cases, too.\\

We lastly show that the case where $\widetilde D_\pm(z_0)$ intersect for some $x_0\in(-1,1)$ will actually happen:\\
For
$$x_0=x^*:=\frac{-2z_0}{1+z_0^2},$$
(\ref{kaihe}) reads
$$\log(1+r)+\log(1-r)-\log r=\log(1+z_0)+\log(1-z_0)-\log z_0-T:=Y\in\R,$$
which means 
$$r=\frac{\sqrt{4+e^{2Y}}-e^Y}{2}.$$
Since $r\left(T,\frac{-2z_0}{1+z_0^2}\right)$ increases with growing $T$, and $r\left(T,\frac{-2z_0}{1+z_0^2}\right)\to1$ for $T\to\infty$, it will at some point of time $T$ become so large that
$$\arctanh r\left(T,\frac{-2z_0}{1+z_0^2}\right) = \frac{\pi}{2}+\arctanh z_0.$$
Then $\sigma(T,x^*)=2\cdot (\arctanh r(T,x^*)-\arctanh z_0)=\pi$ and consequently the curves $\widetilde D_\pm(z_0)$ intersect on $\R^-.$\\

This concludes the proof of The
orem \ref{th2}.\\

We finally prove Theorem \ref{th3} by applying the maximum principle to equation \eqref{Loewnerinvers0} in the free end time version. We have $$\mathcal W_T(z_0) = \{w(T)\,:\, w:[0,\infty)\to\D \text{ solves \eqref{Loewnerinvers0}}, \, T\in[0,\infty)\}.$$
If $w(t)$ is a solution with $w(T)\in \partial \mathcal W_T(z_0),$ then we have the same setting as above and the additional information that
$$\Re H'(\mu_t, \lambda(t), w(t)) = \max_{\mu \in \mathcal P} \Re H'(\mu, \lambda(t), w(t)) = 0$$
for almost all $t\in [0,T],$ see, e.g., \cite{Lewis}, Theorem 5.18.\\
% cf. Th. II.6, \cite{Roth}.\\

The optimal driving term corresponding to \eqref{OptimalDrive} thus has to fulfill 
$$\cos \phi=-\frac{2|w|}{1+|w|^2},$$
which means

$$x=\frac{-2|w|}{1+|w|^2} ,$$
and thus \eqref{s:l} becomes
\begin{align*}
\frac{d}{dt}|w|=|w|\frac{1-|w|^2}{1+|w|^2},\quad |w(0)|=z_0,
\end{align*}
which is equivalent to
$$|w(t)|=\frac{-1+z_0^2 + \sqrt{(1-z_0^2)^2+4z_0^2e^{2t}}}{2e^tz_0}.
$$
We have
$$\frac{d}{dt}\arg w(t)=\pm \frac{2|w|}{1+|w|^2},$$
which yields
$$\frac{d}{d|w|}\arg w= \pm \frac{2}{1-|w|^2},$$
or
$$\arg w=\pm2\left(\arctanh|w|-\arctanh z_0\right)=\pm(d_\D(0,|w|)-d_\D(0,z_0)).$$
Taking into account our results about the sets $W_T(z_0)$, we conclude that $\mathcal W(z_0)=\overline{\mathcal W(z_0)}\cap \D$ and that $\overline{\mathcal W(z_0)}$ is the closed region bounded by $\partial\D$ and the hyperbolic spirals
$$S_\pm(z_0)=\{re^{\pm i\sigma}\,:\,\sigma=d_\D(0,r)-d_\D(0,z_0), \sigma\in[0,\pi]\}.$$
This concludes the proof.

\def\cprime{$'$}
\providecommand{\bysame}{\leavevmode\hbox to3em{\hrulefill}\thinspace}
\providecommand{\MR}{\relax\ifhmode\unskip\space\fi MR }
\providecommand{\MRhref}[2]{%
  \href{http://www.ams.org/mathscinet-getitem?mr=#1}{#2}
}
\providecommand{\href}[2]{#2}


\begin{thebibliography}{{Rog}34}

\bibitem[Dur83]{Duren:1983}
P.~L. Duren, \emph{Univalent functions}, Springer-Verlag, New York, 1983.

\bibitem[GG76]{MR0453993}
V.~V. Gorja{\u\i}nov and V.~Ja. Gutlyanski{\u\i}, \emph{Extremal problems in
  the class {$S_{M}$}}, Mathematics collection ({R}ussian), Izdat. ``Naukova
  Dumka'', Kiev, 1976, pp.~242--246.

\bibitem[Lew06]{Lewis}
A.Lewis, \emph{The Maximum Principle of Pontryagin
in control and in optimal control}, Lecture notes from a course held in the Department of Applied Mathematics in the Polytechnic University of Catalonia, May 2006,
http://www.mast.queensu.ca/~andrew/teaching/MP-course/pdf/maximum-principle.pdf
(access date: 13 August 2015)

\bibitem[LM86]{MR889459}
E.~B. Lee and L.~Markus, \emph{Foundations of optimal control theory}, second
  ed., Robert E. Krieger Publishing Co., Inc., Melbourne, FL, 1986.

\bibitem[MR05]{MarshallRohde:2005}
D.~E. Marshall and S.~Rohde, \emph{The {L}oewner differential equation and slit
  mappings}, J. Amer. Math. Soc. \textbf{18} (2005), no.~4, 763--778.

\bibitem[Pom75]{Pom:1975}
Ch. Pommerenke, \emph{Univalent functions}, Vandenhoeck \& Ruprecht,
  G\"ottingen, 1975.

\bibitem[PS15]{MR3334955}
Dmitri Prokhorov and Kristina Samsonova, \emph{Value range of solutions to the
  chordal {L}oewner equation}, J. Math. Anal. Appl. \textbf{428} (2015), no.~2,
  910--919.

\bibitem[{Rog}34]{rogo1}
Werner {Rogosinski}, \emph{{Zum Schwarzschen Lemma.}}, {Jahresber. Dtsch.
  Math.-Ver.} \textbf{44} (1934), 258--261.

\bibitem[Rot98]{Roth}
O.~Roth, \emph{Control {T}heory in ${H}(\mathbb{D})$}, Ph.D. thesis, University
  of {W}\"urzburg, 1998.

\bibitem[RS14]{MR3262210}
Oliver Roth and Sebastian Schlei{\ss}inger, \emph{Rogosinski's lemma for
  univalent functions, hyperbolic {A}rchimedean spirals and the {L}oewner
  equation}, Bull. Lond. Math. Soc. \textbf{46} (2014), no.~5, 1099--1109.

\end{thebibliography}
\end{document}